\documentclass{amsart}

\usepackage[margin=1.4in]{geometry}

\usepackage{type1cm}         
\usepackage{graphicx}        
\usepackage{multicol}        
\usepackage[bottom]{footmisc}

\usepackage{amsfonts, amsmath, amssymb}
\usepackage{xspace}          
\usepackage{ytableau}        
\usepackage{pgf}

\usepackage[draft=true]{minted} 
\setminted{breaklines=true}
\definecolor{lbcolor}{rgb}{0.9,0.9,0.9}
\setminted{bgcolor=lbcolor}
\setminted{fontsize=\footnotesize}

\usepackage{mathtools}
\usepackage{booktabs}        
\usepackage{tikz}            
\usepackage{csquotes}        

\usepackage[colorinlistoftodos, bordercolor=orange, backgroundcolor=orange!20, linecolor=orange, textsize=scriptsize]{todonotes} 

\usepackage[noend]{algorithmic}

\usepackage{calc,ifthen,alltt,tabularx,capt-of}
\usepackage[linesnumbered,longend,ruled,vlined]{algorithm2e}

\usepackage[bookmarks=false, hidelinks]{hyperref}

\usepackage{longtable}
\usepackage{makecell}
\usepackage{wrapfig}

\usepackage[capitalise, noabbrev]{cleveref} 
\usepackage{enumitem} 

\newcommand\OSCAR{\texttt{OSCAR}\xspace}

\newcommand\polyDB{\texttt{polyDB}\xspace}

\newcommand\FF{\mathbb{F}}
\newcommand\ZZ{\mathbb{Z}} 
\newcommand\QQ{\mathbb{Q}} 
\newcommand\PP{\mathbb{P}}
\newcommand\RR{\mathbb{R}} 
\newcommand\CC{\mathbb{C}} 

\definecolor{promptColor}{rgb}{0.0,0.0,0.589}
\definecolor{brkpromptColor}{rgb}{0.589,0.0,0.0}
\definecolor{gapinputColor}{rgb}{0.589,0.0,0.0}
\definecolor{gapoutputColor}{rgb}{0.0,0.0,0.0}
\definecolor{darkgreen}{rgb}{0.05,0.6,0.1}
\definecolor{colrem}{rgb}{0,0.7,0}

\newcommand{\tmfloatcontents}{}
\newlength{\tmfloatwidth}
\newcommand{\tmfloat}[5]{
	\renewcommand{\tmfloatcontents}{#4}
	\setlength{\tmfloatwidth}{\widthof{\tmfloatcontents}+1in}
	\ifthenelse{\equal{#2}{small}}
	{\setlength{\tmfloatwidth}{0.45\linewidth}}
	{\setlength{\tmfloatwidth}{\linewidth}}
	\begin{minipage}[#1]{\tmfloatwidth}
		\begin{center}
			\tmfloatcontents
			\captionof{#3}{#5}
		\end{center}
\end{minipage}}

\theoremstyle{theorem}
\newtheorem{theorem}{Theorem}[section]

\theoremstyle{definition}
\newtheorem{definition}[theorem]{Definition}
\newtheorem{example}[theorem]{Example}

\newcommand{\fd}{./}

\title{Matroids in OSCAR}
\author{Daniel Corey}
\address{Department of Mathematical Sciences, University of Nevada, Las Vegas, USA}
\email{daniel.corey@unlv.edu}

\author{Lukas K\"uhne}
\address{Fakult\"at f\"ur Mathematik, Universit\"at Bielefeld, Germany}
\email{lukas.kuehne@math.uni-bielefeld.de}

\author{Benjamin Schr\"oter}
\address{Department of Mathematics, KTH Royal Institute of Technology, Stockholm, Sweden}
\email{schrot@kth.se}

\begin{document}
\begin{abstract}
\OSCAR \cite{OSCAR} is an innovative new computer algebra system which combines and extends the power of its four cornerstone systems - GAP (group theory), Singular (algebra and algebraic geometry), Polymake (polyhedral geometry), and Antic (number theory). Here, we present parts of the module handeling matroids in \OSCAR, which will appear as a chapter of the upcoming \OSCAR book \cite{OSCAR-book}.
A matroid is a fundamental and actively studied object in combinatorics. Matroids generalize linear dependency in vector spaces as well as many aspects of graph theory. Moreover, matroids form a cornerstone of tropical geometry and a deep link between algebraic geometry and combinatorics.
Our focus lies in particular on computing the realization space and the Chow ring of a matroid.
\end{abstract}

\maketitle

\section{Introduction}
\noindent A matroid is a combinatorial abstraction of independence, e.g., linear independence of vectors or spanning sets of edges in a graph. Matroids were introduced by Whitney~\cite{Whi35} and independently by Nakasawa (see \cite{Nak09}).  They are a central object in mathematics connecting multiple disciplines such as combinatorics, algebra, and geometry. Consequently, matroids have found applications in fields such as control theory, optimization, and algebraic geometry, see for example~\cite{Bra11,Edm03,GGMS87} respectively. Of particular importance is the work of  Adiprasito, Huh and Katz~\cite{AHK}, who demonstrate that matroids admit a Hodge theory originally stemming from algebraic geometry.

In this chapter, we demonstrate how one can use our implementation of matroids in \OSCAR by focusing on realization spaces and Chow rings of matroids.

\section{Basics}\noindent
Let $V$ be a $d$-dimensional vector space, $E$ a finite set, and $(v_i)_{i\in E}$ a sequence of spanning vectors of $V$. There are various ways to record the linear dependencies among these vectors.

\begin{itemize}
	\item \textit{Independent sets}: $\mathcal{I} = \{A\subseteq E \, : \, (v_i)_{i\in A} \text{ are linearly independent}\}$.
	\item \textit{Bases}: $\mathcal{B} = \{A\subseteq E \, : \, (v_i)_{i\in A} \text{ is a basis of } V\}$.
	\item \textit{Rank}: $\mathrm{rk}: 2^E \to \ZZ_{\geq 0}$; $\mathrm{rk}(A) = \dim \mathsf{span}(v_i \, : \, i\in A)$.
	\item \textit{Flats}: $\mathcal{F} = \{ A \subseteq E \, : \, v_{j} \notin \mathsf{span}(v_i \, : \, i\in A) \text{ for all } j\notin A \}$.
	\item \textit{Circuits}: $\mathcal{C} = \{A\subseteq E \, : \,  (v_i)_{i\in A} \text{ are minimally linearly dependent}  \}$.
\end{itemize}

The notions of independent sets,  bases, rank, flats, and circuits may each be axiomitized, each of which leads to a definition of a matroid. We favor the description in terms of bases.
\begin{definition}
	A \textit{matroid} $\mathsf{M}$ consists of a finite set $E$ and a non-empty collection $\mathcal{B} \subset 2^{E}$ that satisfies the \textit{basis exchange axiom}: for each pair $A,B$ of distinct elements of $\mathcal{B}$ and $x\in A\setminus B$, there is a $y\in B\setminus A$ such that $A\setminus \{x\} \cup \{y\}$ is in $\mathcal{B}$.
\end{definition}
Elements of $\mathcal{B}$ are called \textit{bases} of $\mathsf{M}$. To emphasize dependence on $\mathsf{M}$, especially when $\mathsf{M}$ is defined using a different set of axioms, write $\mathcal{B}(\mathsf{M})$ for $\mathcal{B}$. All bases of $\mathsf{M}$ have the same size \cite[Lemma 1.2.1]{Oxl11}; this common size is called the \emph{rank} of $\mathsf{M}$. The remaining ways of recording dependence can be recovered from this data.

\begin{itemize}
	\item \textit{Independent sets}: $\mathcal{I}(\mathsf{M}) = \{A\subseteq \{1,2,\ldots,n\} \, : \, A\subseteq B \text{ for some } B\in \mathcal{B}\}$.
	\item \textit{Rank}: $\mathrm{rk}_{\mathsf{M}}: 2^E \to \ZZ_{\geq 0}$, $\mathrm{rk}(A) = \max(|I\cap A|\; : \; I\in \mathcal{I})$.
	\item \textit{Flats}: $\mathcal{F}(\mathsf{M}) = \{ A \subseteq E \, : \, \mathrm{rk}(A\cup e) > \mathrm{rk}(A) \text{ for all } e\in E\setminus A \}$.
	\item \textit{Circuits}: $\mathcal{C}(\mathsf{M}) = \{A\subseteq E \, : \,  A\not\in \mathcal{I}(\mathsf{M})  \text{ and } A'\in \mathcal{I}(\mathsf{M})  \text{ for all } A'\subsetneq A\}$.
\end{itemize}
 
For example, consider the matroid $\mathsf{M}$ on $E=\{1,2,3,4\}$ whose bases are $\mathcal{B}(\mathsf{M}) = \{12, 13, 14, 23, 24\}$ where we used, e.g., $12$ as a short hand notation for the set $\{1,2\}$. 
\begin{minted}{jlcon}
julia> B = [[1,2],[1,3],[1,4],[2,3],[2,4]];

julia> M = matroid_from_bases(B,4)
Matroid of rank 2 on 4 elements
\end{minted}
Once created, we can determine other axiomatic characterizations of $\mathsf{M}$.
\begin{minted}{jlcon}
julia> independent_sets(M)
10-element Vector{Vector{Int64}}:
 []
 [1]
 [2]
 [3]
 [4]
 [1, 4]
 [2, 4]
 [1, 3]
 [2, 3]
 [1, 2]

julia> flats(M)
5-element Vector{Vector{Int64}}:
 []
 [1]
 [2]
 [3, 4]
 [1, 2, 3, 4]

julia> circuits(M)
3-element Vector{Vector{Int64}}:
 [3, 4]
 [1, 2, 3]
 [1, 2, 4]
\end{minted}

\begin{figure}
	\includegraphics[width=\textwidth]{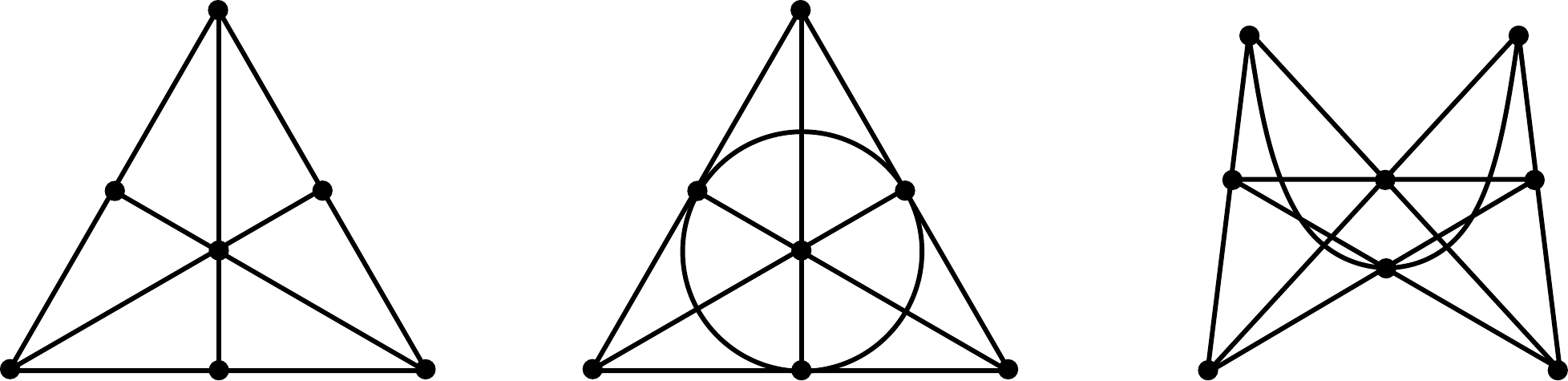}
	\caption{
	Projective realizations of the Fano, non-Fano and M\"{o}bius--Kantor matroids. Points connected by an arc are interpreted as colinear.}
	\label{fig:matroids}
\end{figure}

Consider a vector configuration whose elements are the columns of a $r\times n$ full-rank matrix $X$. The matroid of this configuration, denoted $\mathsf{M}[X]$, is the matroid whose ground set is $\{1,2,\ldots,n\}$ and its bases are the collections of $r$ columns that are of full rank. Because multiplying a column of $X$ by a nonzero number does not change the matroid $\mathsf{M}[X]$, we may also view the columns as points in the projective space $\mathbb{P}^{r-1}$, and so $X$ may be viewed as a projective realization of $\mathsf{M}$. This perspective is particularly useful for illustrating matroids, and projective realizations of some important matroids are depicted in Figure \ref{fig:matroids}. For example, the famous  \textit{Fano matroid} $\mathsf{F}$ is the left-most diagram in this figure. This matroid is $\mathsf{F} = \mathsf{M}[X]$ where $X$ is the $3\times 7$ matrix whose columns are the $7$ points of $\PP^2(\FF_2)$. 
\begin{minted}{jlcon}
julia> X = matrix(GF(2), [1 0 1 0 1 0 1; 0 1 1 0 0 1 1; 0 0 0 1 1 1 1])
[1   0   1   0   1   0   1]
[0   1   1   0   0   1   1]
[0   0   0   1   1   1   1]


julia> F = matroid_from_matrix_columns(X)
Matroid of rank 3 on 7 elements
\end{minted}
Let $\mathsf{M}_1$ and $\mathsf{M}_2$ be matroids with ground sets $E_{1}$ and $E_{2}$, respectively. An \textit{isomorphism} of matroids $\varphi:\mathsf{M}_1 \to \mathsf{M}_2$ is a bijection $\varphi:E_1 \to E_2$ such that $A\in \mathcal{B}(\mathsf{M}_1)$ if and only if $\varphi(A) \in \mathcal{B}(\mathsf{M}_2)$. We may verify that $\mathsf{M}[X]$ is  isomorphic to the Fano matroid. 
\begin{minted}{jlcon}
julia> is_isomorphic(F, fano_matroid())
true
\end{minted}
The \textit{automorphism group} of a matroid $\mathsf{M}$, denoted $\mathsf{Aut}(\mathsf{M})$, is the group of all isomorphism from $\mathsf{M}$ to itself. The automorphism group of the Fano matroid is $\mathsf{PSL}_{3}(\FF_2)$.
\begin{minted}{jlcon}
julia> G = automorphism_group(F);

julia> describe(G)
"PSL(3,2)"
\end{minted}

\section{Realizability}\noindent 
In the previous section, we saw how to obtain a matroid from a vector configuration, and when the vectors are arranged as columns of a matrix $X$, the resulting matroid is $\mathsf{M}[X]$. Let $\FF$ be a field. A matroid $\mathsf{M}$ is $\FF$-\textit{realizable} if there is a matrix~$X$ with entries in the field $\FF$ such that $\mathsf{M} \cong \mathsf{M}[X]$. If $\mathsf{M}$ is $\FF$-realizable for some field $\FF$, then $\mathsf{M}$ is said to be realizable. 

Denote by $\binom{\{1,\ldots,n\}}{r}$ the $r$-element subsets of $\{1,\ldots,n\}$. The smallest matroid, with respect to the size of the ground set,
that is not realizable over any field is the \textit{V\'amos} matroid.
\begin{minted}{jlcon}
julia> nonbases(vamos_matroid())
5-element Vector{Vector{Int64}}:
 [1, 2, 3, 4]
 [1, 2, 5, 6]
 [1, 2, 7, 8]
 [3, 4, 5, 6]
 [5, 6, 7, 8]
\end{minted}
It is a rank $4$ matroid on a set of size  $8$, and its bases (up to isomorphism) are the elements of $\binom{\{1,\ldots,8\}}{4}$ other than those listed above. We verify below with our code that it is indeed not realizable.

Determining whether a matroid is realizable is in general a hard problem. There does not exist a finite list of axioms that combinatorially characterize realizability \cite{MayhewNewmanWhittle}, which led V\'amos to quip ``the missing axiom of matroid theory is lost forever'' \cite{Vamos}.  Nevertheless, the V\'amos matroid violates Ingleton's inequality \cite{Ing71} and hence is not realizable; see \cite[\S 6.1, Exer. 7]{Oxl11}.
Moreover, in \cite{Sturmfels:1987}, Sturmfels proves that matroid-realizability over $\QQ$ is equivalent to the solvability of Diophantine equations over $\QQ$. This is  an extension of Hilbert's 10th problem to the rationals and it is not known whether this problem is decidable or not.

On the other hand, as a result of the discussion below, realizability of matroids over algebraically closed fields can be decided using Gr\"obner bases. Nevertheless, there are practical limitations in that this algorithm becomes rather slow if the matroid is relatively large.

A closely related concept to realizability is the notion of a realization space of the matroid~$\mathsf{M}$. Given a field $\FF$ and a matroid $\mathsf{M}$, its \textit{realization space}~$\mathcal{R}(\mathsf{M};\FF)$ is a (possibly empty) algebraic variety (more generally, a scheme), defined over $\FF$, whose closed points parameterize equivalence classes of point configurations in $\mathbb{P}_{\FF}^{r-1}$ whose matroid is $\mathsf{M}$, where two configurations are equivalent if one can be transformed to the other by an element of $\mathsf{PGL}_r(\FF)$. In particular, the matroid $\mathsf{M}$ is $\FF$-realizable if and only if~$\mathcal{R}(\mathsf{M};\FF) \neq \emptyset$.
\subsection{Algorithmic Framework}

In this subsection we describe the basic algorithmic framework to compute the affine coordinate ring of a matroid realization space. We illustrate this code in the a series of examples in the next subsection.  Following~\cite[Thm. 6.8.9]{Oxl11} the basic algorithm does the following for a rank $r$ matroid $\mathsf{M}$ on the ground set~$\{1,\ldots,n\}$. Assume that $\mathsf{M}$ has no loops (a \textit{loop} of a matroid is an element not contained in any basis). To streamline the exposition, we consider the special case where $B_{0} = \{n-r+1,\ldots,n\}$ is a basis of $\mathsf{M}$. 
Let $A = \ZZ[x_{ij}]$ be the polynomial ring in the $r(n-r)$ determinants $x_{ij}$ for $1\leq i \leq r$, $1 \leq j \leq n-r$. Define the matrix
\begin{equation*}
X = \begin{bmatrix}
x_{11} & x_{12} & \cdots & x_{1,n-r} & 1 & 0 & \cdots & 0  \\
x_{21} & x_{22} & \cdots & x_{2,n-r} & 0 & 1 & \cdots & 0  \\
\cdots & \cdots & \cdots & \cdots & \cdots & \cdots & \cdots & \cdots  \\
x_{r1} & x_{r2} & \cdots & x_{r,n-r} & 0 & 0 & \cdots & 1  
\end{bmatrix}.
\end{equation*}
The matroid $\mathsf{M}$ is realizable over a field $\FF$ if we can replace every $x_{ij}$ with an element of $\FF$ such that $\mathsf{M}=\mathsf{M}\left[X\right]$. For each $r$-element subset $N$ of~$\{1,\ldots,n\}$, write $p_{N}$ for the determinant of $X$ whose columns are indexed by the elements of $N$ (preserving the order). To realize $\mathsf{M}$, the $x_{ij}$ must satisfy $p_{N} = 0$ for $N\in \binom{\{1,\ldots,n\}}{r} \setminus \mathcal{B}(\mathsf{M})$ and $p_{N} \neq 0$ for $N \in \mathcal{B}(\mathsf{M})$.  Observe that $x_{ij} = (-1)^{i+1}p_{B}$ for $i\in \{1,\ldots,r\}$ and $j\in \{1,\ldots,n-r\}$ where $B = B_{0} \setminus \{n-r+i\} \cup \{j\}$. Therefore, we must have $x_{ij} =0$ whenever $B_{0} \setminus \{n-r+i\} \cup \{j\}\not\in \mathcal{B}(\mathsf{M})$. 

The matroid $\mathsf{M}[X]$ is invariant under row and column scaling by a nonzero scalar as this does not change the linear dependence structure among the columns. In practice, this means that we may set equal to 1 the first nonzero entry of each column, and then among the remaining entries, we may set equal to 1 the first nonzero entry of each row (if this entry exists).

Let us carefully carry out this procedure. Let us denote the row with the first nonzero entry in column $j\in \{1,\ldots,n-r\}$ by $\mu(j)$ where $\mu(j)=r+1$ if there is no such entry. In other words we obtain the function $\mu:\{1,\ldots,n-r\} \to \{1,\ldots,r+1\}$ with
\begin{equation*}
\mu(j) = \mathrm{min}(\{r+1\} \cup \{i\in \{1,\ldots,r\} \, : \, B_{0} \setminus \{n-r+i\} \cup \{j\}\in \mathcal{B}(\mathsf{M})\}).
\end{equation*}
Similarly the first nonzero (column) entry $j$ in row $i\in \{1,\ldots,r\}$ not fulfilling $\mu(j)=i$ is denoted by $\nu(i)$. That is
$\nu:\{1,\ldots,r\} \to \{1,\ldots,n-r+1\}$ is the function
\begin{equation*}
\nu(i) = \mathrm{min}( \{n-r+1\} \cup \{j\in \{1,\ldots,n-r\} \, : \, B_{0} \setminus \{n-r+i\} \cup \{j\}\in \mathcal{B}(\mathsf{M}) \text{ and }\mu(j)\neq i\}).
\end{equation*}
where $\nu(i)=n-r+1$ indicates that there is no such entry.
Let us further set
\begin{equation*} 
A_{\mathsf{M}} = \ZZ[x_{ij}\, : \, B_{0} \setminus \{n-r+i\} \cup \{j\}\in \mathcal{B}(\mathsf{M}), i\neq \mu(j), j\neq \nu(i)]
\end{equation*}
Define a ring homomorphism
\begin{equation*}
\pi: A \to A_{\mathsf{M}} \hspace{10pt} \pi(x_{ij}) = 
\begin{cases}
0 & \text{if } B_{0} \setminus \{n-r+i\} \cup \{j\}\not \in \mathcal{B}(\mathsf{M}), \\
1 & \text{if }  i = \mu(j) \text{ or } j = \nu(i), \\
x_{ij}  & \text{otherwise}.
\end{cases}
\end{equation*}
Define the ideal $I_{\mathsf{M}}$ and multiplicative semigroup $U_{\mathsf{M}}$ by
\begin{equation*}
I_{\mathsf{M}} = \langle \pi(p_{N})  \mid   N\notin \mathcal{B}(\mathsf{M}), \,  |N|=r \rangle, \text{ and } U_{\mathsf{M}} = \langle \pi(p_{B})  \mid  B \in \mathcal{B}(\mathsf{M}) \rangle_{\mathrm{smgp}}.
\end{equation*}
Here, the subscript $\mathrm{smgp}$ indicates that we take the multiplicative semigroup generated by the listed elements.
The realization space $\mathcal{R}(\mathsf{M};\mathbb{F})$ is an affine scheme, and its affine coordinate ring is $S_{\mathsf{M}} \otimes_{\ZZ} \mathbb{F}$ where
\begin{equation*} 
S_{\mathsf{M}} = U_{\mathsf{M}}^{-1} A_{\mathsf{M}} / I_{\mathsf{M}}.
\end{equation*}
The matroid $\mathsf{M}$ is realizable if and only if $\mathcal{R}(\mathsf{M}, \mathbb{F})$ is nonempty, equivalently,  $S_{\mathsf{M}}$ is not the zero ring. In practice, showing $S_{\mathsf{M}}$ is the zero ring amounts to showing that the saturation of $I_{\mathsf{M}}$ with respect to the semigroup $S_{\mathsf{M}}$ is the unit ideal in $A_{\mathsf{M}}$.  When~$\FF$ is algebraically closed, a closed point in $\mathcal{R}(\mathsf{M};\FF)$ corresponds to a realization of $\mathsf{M}$ via the matrix $X$.

The details of this implementation are also described in~\cite{BK21}. To speed up this computation we use the following techniques.
\begin{itemize}
	\item The assumption that  $B_0=\{n-r+1,\dots,n\}$ is a basis of $\mathsf{M}$ allows us to assume that the last $r$ columns of $X$ form the $r\times r$ identity matrix. Of course, the above procedure may be modified so that the submatrix formed by the columns of $B_{0}$ is the identity matrix \textit{for any} $B_{0} \in \mathcal{B}(\mathsf{M})$.
	\item
	We do a short check before the realization space computation for which bases the ring $A_{\mathsf{M}}$ has the fewest variables and use this basis instead as this makes the following computations easier.
	\item We embed the space $\mathcal{R}(\mathsf{M};\mathbb{F})$ into a smaller dimensional ambient space using techniques described in~\cite[Appendix~A]{BK21} and \cite[\S~6]{CL22}.
	\item We simplify the ideal $I_{\mathsf{M}}$ by computing a reduced Gr\"obner basis and reduce the polynomials in $U_\mathsf{M}$ with respect to this basis.
\end{itemize}

\subsection{Examples}\noindent
We illustrate the procedure from the previous section in the following examples. 

\begin{example}
	Let's revisit the Fano matroid $\mathsf{F}$ introduced in the previous section:
\begin{minted}{jlcon}
julia> realization_space(F, simplify=false)
The realization space is
[0   1    1    1    1   0   0]
[1   0    1   x2    0   1   0]
[1   0   x1    0   x3   0   1]
in the Multivariate polynomial ring in 3 variables over ZZ
within the vanishing set of the ideal
ideal(2, x3 + 1, x2 - 1, x1 + 1)
avoiding the zero loci of the polynomials
RingElem[x2, x3, x1, x1*x2 - x2*x3 + x3, x1 - x3 - 1, x1 + x2 - 1]
\end{minted}
	The ring $A_{\mathsf{F}}$ is the polynomial ring in 3 variables with \textit{integer} coefficients. The output shows the matrix $\pi(X)$ which parametrizes the realizations, the ideal $I_\mathsf{F}$ and the inequations in $U_{\mathsf{F}}$. The equations actually already determine the variables $x_1,x_2,x_3$. So we can simplify the realization space using the equations in $I_{\mathsf{F}}$. We use exactly equations of this kind to embed $\mathcal{R}(\mathsf{F};\FF)$ into a smaller dimensional ambient space as mentioned above.
\begin{minted}{jlcon}
julia> realization_space(F, simplify=true)
The realization space is
[0   1   1   1   1   0   0]
[1   0   1   1   0   1   0]
[1   0   1   0   1   0   1]
in the Integer ring
within the vanishing set of the ideal
2ZZ
\end{minted}
	Because the vanishing ideal is $\langle 2 \rangle$, the matroid is realizable only over fields of characteristic $2$.
	Over any such field a realization of the Fano matroid is given by this matrix.
	We can also check directly that this matroid is not realizable over fields of other characteristics. For instance, we may check characteristic~$5$:
\begin{minted}{jlcon}
julia> is_realizable(F, char=5)
false
\end{minted}
\end{example}
The default option for the parameter \texttt{simplify} is \texttt{true} so we will omit the parameter from now on.

\begin{example}
	A matroid closely related to the Fano matroid $\mathsf{F}$ is the non-Fano matroid. This matroid is obtained from the $\mathsf{F}$ by turning one of the seven circuits of $\mathsf{F}$ into a basis, see Figure~\ref{fig:matroids}.
\begin{minted}{jlcon}
julia> realization_space(non_fano_matroid())
The realization space is
[1   1   0   0   1   1   0]
[0   1   1   1   1   0   0]
[0   1   1   0   0   1   1]
in the Integer ring
avoiding the zero loci of the polynomials
RingElem[2]
\end{minted}
	In this case, the defining ideal is trivial, but the element $2$ must not vanish in the field we realize this matroid. Therefore the matroid is realizable over fields of every characteristic except $2$.
\end{example}

\begin{example}
	We can also confirm that the V\'amos matroid is not realizable over any field:
\begin{minted}{jlcon}
julia> is_realizable(vamos_matroid())
false
\end{minted}
\end{example}

\begin{example}
	Let $\mathsf{M}$ be the \textit{M\"{o}bius-Kantor} matroid. 
\begin{minted}{jlcon}
julia> M = moebius_kantor_matroid();

julia> nonbases(M)
8-element Vector{Vector{Int64}}:
 [1, 2, 3]
 [1, 4, 5]
 [1, 7, 8]
 [2, 4, 6]
 [2, 5, 8]
 [3, 4, 7]
 [3, 6, 8]
 [5, 6, 7]
\end{minted}
	This is a rank-$3$ matroid on $\{1,\ldots,8\}$ whose bases (up to isomorphism) are the elements of $\binom{\{1,\ldots,8\}}{3}$ other than those listed above, see Figure \ref{fig:matroids} for an illustration. We may compute a simple presentation of the ring $S_{\mathsf{M}}$ which defines an embedding of $\mathcal{R}(\mathsf{M};\FF)$ in a smaller dimensional space:
\begin{minted}{jlcon}
julia> realization_space(M)
The realization space is
[1   0   1   0   1     0         1    1]
[0   1   1   0   0     1         1   x1]
[0   0   0   1   1   -x1   -x1 + 1    1]
in the Multivariate polynomial ring in 1 variable over ZZ
within the vanishing set of the ideal
ideal(x1^2 - x1 + 1)
avoiding the zero loci of the polynomials
RingElem[x1, x1 - 1]
\end{minted}
	This means that 
	\begin{equation*} 
	S_{\mathsf{M}} \cong \ZZ[x^{\pm}] / \langle x^2-x+1 \rangle.
	\end{equation*}
	(The polynomial $x-1$ is already a unit in this ring.)
	Therefore, $\mathsf{M}$ is $\FF$-realizable if and only if $x^2-x+1$ has a root in $\FF$ (equivalently, $\FF$ has a primitive 6-th root of unity). In characteristic 0, we see that $\mathsf{M}$ is not realizable over $\RR$, and therefore not orientable.
	
	We can also check whether $\mathsf{M}$ is realizable over a specific finite field $\mathbb{F}_q$ for some prime power $q$.
\begin{minted}{jlcon}
julia> is_realizable(M, q=9)
true
\end{minted}
	Hence $\mathsf{M}$ is realizable over $\mathbb{F}_9$.
	We can also check this systematically for all prime powers of size up to  $13$.
\begin{minted}{jlcon}
julia> qs = [2,3,4,5,7,8,9,11,13];
julia> println([is_realizable(M, q=a) for a in qs])
Bool[0, 1, 1, 0, 1, 0, 1, 0, 1]
\end{minted}
	From this list we can deduce the hypothesis that $\mathsf{M}$ is realizable over the finite field $\mathbb{F}_q$ if and only if $q\not \equiv 2 \mod 3$.
\end{example}

As mentioned in the introduction, a second motivation for matroids stems from graph theory. An undirected graph $G=(V,E)$ gives rise to a \emph{graphic matroid} $\mathsf{M}(G)$ on the ground set $E$ and bases the maximal subsets of edges that don't contain a cycle.
\begin{example}
In this example, we compute the realization space of the graphic matroid of the complete graph on four vertices $K_4$. 
\begin{minted}{jlcon}
julia> M = cycle_matroid(Graphs.complete_graph(4));
julia> realization_space(M)
The realizations are parametrized by
[1   0   1   0   1    0]
[0   1   1   0   0    1]
[0   0   0   1   1   -1]
in the Integer Ring
\end{minted}
	Hence there are no algebraic conditions on the field, and this matrix defines a realization of $\mathsf{M}(K_4)$ over every field (of course, $-1=1$ over fields of characteristic 2). That $\mathcal{R}(\mathsf{M}(K_4);\FF)$ consists of a single point for any field $\FF$ is consistent with the fact due to White which asserts that $\mathcal{R}(\mathsf{M};\FF)$ consists of a single point whenever $\mathsf{M}$ is realizable over $\FF_2$, see \cite[Theorem~6.3]{Kat16}.
\end{example}

\begin{example}
	Another prominent matroid is the \emph{Pappus} matroid as its configurations of rank two flats is a Pappus configuration. Let's begin by computing its realization space over $\mathbb{C}$:
\begin{minted}{jlcon}
julia> RS = realization_space(pappus_matroid(), char=0)
The realization space is
[1   0   1   0   x2   x2                 x2^2    1    0]
[0   1   1   0    1    1   -x1*x2 + x1 + x2^2    1    1]
[0   0   0   1   x2   x1                x1*x2   x1   x2]
in the Multivariate polynomial ring in 2 variables over QQ
avoiding the zero loci of the polynomials
RingElem[x1 - x2, x2, x1, x2 - 1, x1 + x2^2 - x2,
    x1 - 1, x1*x2 - x1 - x2^2]
\end{minted}
	Thus the realization space is the affine space $\mathbb{A}^2$ over $\mathbb{C}$ with the seven specified curves removed.
	One can obtain a specific realization of this matroid by picking a point in that space that avoids these seven exceptional curves:
\begin{minted}{jlcon}
julia> realization(RS)
One realization is given by
[1   0   1   0   2   2   4   1   0]
[0   1   1   0   1   1   1   1   1]
[0   0   0   1   2   3   6   3   2]
in the Rational field
\end{minted}
\end{example}

\begin{example}
	\OSCAR is connected to the database \polyDB\cite{polyDB}. To illustrate this functionality, we collect all simple matroids with rank $3$ on $9$ elements, up to isomorphism.
        \footnote{The database contains one representative for every isomorphism class of matroids.}
\begin{minted}{jlcon}
julia> db = Polymake.Polydb.get_db();
julia> collection = db["Matroids.Small"];
julia> query = Dict("RANK"=>3,"N_ELEMENTS"=>9,"SIMPLE"=>true);
julia> results = Polymake.Polydb.find(collection, query);
julia> oscar_matroids = [Matroid(pm) for pm in results];
julia> length(oscar_matroids)
383
\end{minted}
	There are $383$ such matroids, and only $370$ of these are realizable over a field of characteristic $0$.
\begin{minted}{jlcon}
julia> char_0_matroids = [M for M in oscar_matroids 
        if is_realizable(M, char=0)];
julia> length(char_0_matroids)
370
\end{minted}
\end{example}

As a consequence of Mn\"ev's universality theorem \cite{Mne85, Mne88}, realization spaces for rank-$3$ matroids satisfies Murphy's law in the sense of Vakil \cite{Vak06}, i.e., each singularity type appears in the realization space of a rank-$3$ matroids contains. While the modern proofs of this fact are constructive \cite{Car15, Laf03, LV13}, the ground sets of the matroids produced are large, even for the simplest singularities. In \cite{CL23}, the authors use these \OSCAR realization space functions to prove that realization spaces for $\CC$--realizable, rank-$3$ matroids on ground sets with fewer than 11 elements are all smooth, but there are rank-$3$ matroids on at least $12$ elements whose realization spaces have nodal singularities.

\section{Chow rings of Geometries} \noindent
To a matroid one might associate an element of an abelian group. Such a mapping is called a \emph{matroid invariant} if it is constant for all members in an matroid isomorphism class (analogous definitions can be made for graphs).  Prominent examples include the Tutte polynomial or the characteristic polynomial. For an overview on matroid invariants and the Tutte polynomial see, e.g., \cite{FS22,Tutte}.

The \emph{Tutte polynomial} of a matroid $\mathsf{M}$  is the bivariate polynomial
\[
T_\mathsf{M}(x,y) = \sum_{A\subseteq E} (x-1)^{\mathrm{rk}(E)-\mathrm{rk}(A)}(y-1)^{|A|-\mathrm{rk}(A)},  
\]
and the \emph{characteristic polynomial} of $\mathsf{M}$ is the specialization
\[
\chi_\mathsf{M}(q) =  (-1)^{\mathrm{rk}(E)}T_{\mathsf{M}}(1-q,0) = \sum_{A\subseteq E} (-1)^{|A|}\, q^{\mathrm{rk}(E)-\mathrm{rk}(A)}.
\]
If $\mathsf{M}$ is the graphic matroid of a simple graph $\mathsf{G}$ with $c$ connected components then $T_\mathsf{M}(x,y)$ is the Tutte polynomial of $\mathsf{G}$ and $q^c \chi_\mathsf{M}(q)$ coincides with the \emph{chromatic polynomial} $\chi_G(q)$ of $\mathsf{G}$ which counts the number of node colorings using $q$ colors such that adjacent nodes have distinct colors.

For example the complete graph $\mathsf{K}_n$ on $n$ nodes has
\[
\chi_{\mathsf{K}_n}(q) = q\cdot (q-1)\cdots (q-n+1)
\]
$q$-colorings. The following lines of code reproduce this fact in \OSCAR for the complete graph on $n=4$ nodes. They show how one can compute the Tutte and chromatic polynomial of the matroid $\mathsf{M}(\mathsf{K}_n)$.
\begin{minted}{jlcon}
julia> M = cycle_matroid(Graphs.complete_graph(4))
Matroid of rank 3 on 6 elements

julia> tutte_polynomial(M)
x^3 + 3*x^2 + 4*x*y + 2*x + y^3 + 3*y^2 + 2*y

julia> char_poly = characteristic_polynomial(M)
q^3 - 6*q^2 + 11*q - 6

julia> factor(char_poly)
1 * (q - 2) * (q - 1) * (q - 3)
\end{minted}

Observe that in this example the coefficients $w_i$ of the characteristic polynomial $\chi_\mathsf{M}(q) = \sum_{j=0}^\mathrm{rk(E)} \omega_j q^j$ form a \emph{log-concave} sequence, i.e., for all $0 < j < \mathrm{rk}(\mathsf{M})$ holds $\omega^2_j \geq \omega_{j-1}\omega_{j+1}$. 
That this holds true for all matroids has been conjectured by Welsh \cite{Welsh76} and was known as the Heron-Rota-Welsh Conjecture. This conjecture was famously proved by Adiprasito, Huh and Katz in \cite{AHK}, which also contains historical background and further details about this conjecture.

To a loopfree matroid $\mathsf{M}$, define the algebra $\mathrm{A}(\mathsf{M})$ by
\[
	\mathrm{A}(\mathsf{M}) = \mathbb{Q}[x_F\,|\, F \text{ is a proper flat of } \mathsf{M}]/(I+J)
\]
where $I$ and $J$ are the ideals
\begin{equation*}
	I = \left\langle x_F\cdot x_G \,|\, F,G \text{ are incomparable}\right\rangle \text{ and } J=\left\langle\sum_{F\ni i} x_F - \sum_{F\ni j} x_F \,|\, i,j\in E \right\rangle.
\end{equation*}
The algebra $\mathrm{A}(\mathsf{M})$ is the \emph{Chow ring} of the matroid $\mathsf{M}$. We obtain it in \OSCAR via the following command.
\begin{minted}{jlcon}
julia> A = chow_ring(M);

julia> GR, _ = graded_polynomial_ring(QQ,symbols(base_ring(A)));

julia> AA = chow_ring(M,ring=GR);
\end{minted}

This ring is graded by the polynomial degrees and generated in degree $1$. In the code the ring \texttt{AA} is the graded version of the ring \texttt{A}. 
The component $\mathrm{A}^{\mathrm{rk}(E)-1}(\mathsf{M})$ is a one-dimensional $\mathbb{Q}$ vector space.
The \emph{volume map} is an isomorphism 
\[
	\mathrm{vol}_{\mathsf{M}}: \mathrm{A}^{\mathrm{rk}(E)-1}(\mathsf{M}) \to \mathbb{Q}
\]
normalized such that $\mathrm{vol}_{\mathsf{M}}( x_{F_1} \cdots x_{F_{\mathrm{rk}(E)-1}} ) = 1$ for one (and therefore all) chains of flats $\emptyset \subsetneq F_1 \subsetneq \dots \subsetneq F_{\mathrm{rk}(E)-1} \subsetneq E$.

We obtain the volume map in \OSCAR via the following call.
\begin{minted}{jlcon}
julia> vol_map = volume_map(M,AA)
(::Oscar.var"#volume_map#2171") (generic function with 1 method)
\end{minted}

The names and definitions are motivated from algebraic and toric geometry. The Chow ring of a $\mathbb{C}$-realizable matroid is the usual Chow ring of the ``wonderful compactification'' of the complement of a hyperplane arrangement realizing $\mathbb{C}$. See the work of de Concini and Procesi \cite{DP95}.

A degree 1 element $\sum c_F x_F \in\mathrm{A}^1(\mathsf{M})$ is called a \emph{Lefschetz element} whenever
\[
	c_F + c_G > c_{F\cap G} + c_{F\cup G} \text{ for all flats $F, G$ of the matroid } \mathsf{M},
\]
with $c_\emptyset = c_E = 0$. For example the elements
\begin{align*}
	\alpha_i = \sum_{F\ni i} x_F \qquad\text{ and }\qquad
	\beta_i = \sum_{F \not\ni i} x_F.
\end{align*}
are limit points of Lefschetz elements as they weakly satisfy the above inequalities.
Notice that neither $\alpha_i$ nor $\beta_i$ depend on the choice of $i\in E$ as $\beta_j-\beta_i=\alpha_i-\alpha_j$, which is contained in $J$.
Thus from now on we denote the elements $\alpha_i$ and $\beta_i$ by $\alpha$ and $\beta$ respectively.

To construct these two elements in \OSCAR we pick all flats that either contain a certain element or exclude that element.
\begin{minted}{jlcon}
julia> e = matroid_groundset(M)[1];

julia> proper_flats = flats(M)[2:length(flats(M))-1];

julia> a = sum([AA[i] for i in 1:length(proper_flats) if e in proper_flats[i]])
x_{1} + x_{1,2,3} + x_{1,4,5} + x_{1,6}

julia> b = sum([AA[i] for i in 1:length(proper_flats) if !(e in proper_flats[i])])
x_{2} + x_{3} + x_{4} + x_{5} + x_{6} + x_{2,5} + x_{2,4,6} + x_{3,4} + x_{3,5,6}
\end{minted}

A Lefschetz element gives rise to the following three properties known as the \emph{K\"ahler package}. 
\begin{theorem}[{\cite[Theorem 1.4 and 6.19]{AHK}}] Let $\ell\in\mathrm{A}^1(\mathsf{M})$ be a Lefschetz element. Then the following three properties hold.

\item[Poincar\'e duality:] For every non-negative integer $k\leq \tfrac{\mathrm{rk}(E)-1}{2}$, the bilinear pairing 
\[
	\mathrm{A}^k(\mathsf{M}) \times \mathrm{A}^{\mathrm{rk}(E)-k-1} (\mathsf{M}) \to \mathbb{Q},\qquad 
	(\eta_1,\,\eta_2) \mapsto \mathrm{vol}_\mathsf{M}( \eta_1\cdot\eta_2  )
\] 
is non-degenerate.

\item[Hard Lefschetz property:]  For every non-negative integer $k\leq \tfrac{\mathrm{rk}(E)-1}{2}$, the multiplication map
\[
        \mathrm{A}^k(\mathsf{M}) \to \mathrm{A}^{\mathrm{rk}(E)-k-1} (\mathsf{M}),\qquad 
        \eta \mapsto \ell^{\mathrm{rk}(E)-2k-1}\cdot\eta
\]
is an isomorphism.

\item[Hodge--Riemann relations:] For every non-negative integer $k\leq \tfrac{\mathrm{rk}(E)-1}{2}$, the bilinear form
\[
	\mathrm{A}^k(\mathsf{M}) \times \mathrm{A}^{k} (\mathsf{M}) \to \mathbb{Q},\qquad
        (\eta_1,\,\eta_2) \mapsto (-1)^k\,\mathrm{vol}_\mathsf{M}( \eta_1\cdot\ell^{\mathrm{rk}(E)-2k-1}\cdot\eta_2  )
\]
is positive definite on the kernel of the multiplication by $\ell^{\mathrm{rk}(E)-2k}$.
\end{theorem}

We now demonstrate one of many ways to verify these properties for $k=1$ and $\ell=\beta$.  While $\beta$ is only the limit of Lefschetz elements, it still satisfies the conclusions of the above theorem.
First we generate a basis for each of the two vector spaces $\mathrm{A}^k(\mathsf{M})$ and $\mathrm{A}^{\mathrm{rk}(E)-k-1}(\mathsf{M})$.
From these bases we derive the matrix \texttt{Mat1} that represents the bilinear pairing of the Poincar\'e duality. 
\begin{minted}{jlcon}
julia> k = 1
1

julia> R = base_ring(AA);

julia> g = grading_group(R)[1];

julia> PD1, mapPD1 = homogeneous_component(AA,k*g);

julia> basis_PD1 = [mapPD1(x) for x in gens(PD1)];

julia> PD2, mapPD2 = homogeneous_component(AA,(rank(M)-k-1)*g);

julia> basis_PD2 = [mapPD2(x) for x in gens(PD2)];

julia> Mat1 = matrix(QQ,[[vol_map(b1*b2) for b1 in basis_PD1] for b2 in basis_PD2])
[-1    0    0    0    0    0    0    1]
[ 0   -1    0    0    0    0    0    0]
[ 0    0   -1    0    0    0    0    1]
[ 0    0    0   -1    0    0    0    0]
[ 0    0    0    0   -1    0    0    1]
[ 0    0    0    0    0   -1    0    0]
[ 0    0    0    0    0    0   -1    0]
[ 1    0    1    0    1    0    0   -2]
\end{minted}
It is easy to see that the matrix \texttt{Mat1} is square and of full rank, and hence $\mathrm{A}^k(\mathsf{M})\cong\mathrm{A}^{\mathrm{rk}(E)-k-1}(\mathsf{M})$.

Next we verify the hard Lefschetz property by looking at the bilinear form $\mathrm{A}^k\mathsf{M})\times\mathrm{A}^k(\mathsf{M}) \to \QQ$, $(\eta_1,\eta_2)\mapsto \mathrm{vol}_\mathsf{M} (\eta_1 \cdot \ell^{\mathrm{rk}(E)-2k-1}\cdot \eta_2)$ which is non-degenerate only if the hard Lefschetz map is an isomorphism as we already know $\mathrm{A}^k(\mathsf{M})\cong\mathrm{A}^{\mathrm{rk}(E)-k-1}(\mathsf{M})$.
\begin{minted}{jlcon}
julia> Mat2 = matrix(QQ,[[vol_map(b1*b^(rank(M)-2k-1)*b2) for b1 in basis_PD1] for b2 in basis_PD1]);

julia> Mat1 == Mat2
true
\end{minted}
This indeed is the case as the representing matrix \texttt{Mat2} agrees with the full rank matrix \texttt{Mat1}. 

The last property that we inspect is the Hodge--Riemann relations for which we begin by finding a basis of the kernel of the map $\eta\mapsto \ell^{\mathrm{rk}(E)-k}\cdot\eta$.
\begin{minted}{jlcon}
julia> RR, _ = graded_polynomial_ring(QQ,"y_#" => 1:length(basis_PD1));

julia> map = hom(RR,AA,basis_PD1);

julia> K = kernel(hom(RR,AA,[b^(rank(M)-2k)*b for b in basis_PD1]));

julia> basis_HR = [map(h) for h in gens(K) if degree(h).coeff==k*g.coeff]
7-element Vector{MPolyQuoRingElem{MPolyDecRingElem{QQFieldElem, QQMPolyRingElem}}}:
 -x_{6} + 2*x_{1,2,3}
 -x_{6} + 2*x_{1,4,5}
 -x_{6} + 2*x_{1,6}
 -x_{6} + 2*x_{2,5}
 -x_{6} + 2*x_{2,4,6}
 -x_{6} + 2*x_{3,4}
 -x_{6} + 2*x_{3,5,6}
\end{minted}
It remains to check that the bilinear form restricted to this kernel is positive definite.
\begin{minted}{jlcon}
[2   6    4   2    4   2    4]
[4   4   10   4    6   4    6]
[2   2    4   6    4   2    4]
[4   4    6   4   10   4    6]
[2   2    4   2    4   6    4]
[4   4    6   4    6   4   10]

julia> is_positive_definite(matrix(ZZ,[ZZ(i) for i in Mat3]))
true
\end{minted}

Our next goal is to draw conclusions on the coefficients of the characteristic polynomial $\chi_\mathsf{M}(q)$ or more precisely the reduced characteristic polynomial $\overline{\chi}_\mathsf{M}(q)=\chi_\mathsf{M}(q)/(q-1)=\sum \overline{\omega}_j q^j$.
The key is the fact that the $j$-th coefficient $\overline{\omega}_j$ of $\overline{\chi}_\mathsf{M}(q)$ agrees with $(-1)^j \mathrm{vol}_{\mathsf{M}}(\alpha^{\mathrm{rk}(E)-j-1}\beta^j)$ and that multiplication with $\alpha$ corresponds to a truncation of $\mathsf{M}$, that is $\alpha \cdot \mathrm{A}(\mathsf{M})\cong \mathrm{A}(\mathsf{M'})$ where the bases of $\mathsf{M'}$ are the corank-$1$ independent sets of $\mathsf{M}$.
\begin{minted}{jlcon}
julia> reduced_characteristic_polynomial(M)
q^2 - 5*q + 6

julia> [vol_map(a^(rank(M)-j-1)*b^j) for j in range(0,rank(M)-1)]
3-element Vector{QQFieldElem}:
 1
 5
 6
\end{minted}
Thus the Hodge--Riemann relations for the rank $2$ truncation of $\mathsf{M}$ reveal that $\overline\omega_1^2 \geq \overline\omega_0\overline\omega_2$ as this is the determinant of the pairing with basis $\{\alpha,\beta\}$.
Iterating this argument leads then to the fact that the coefficients of $\overline{\chi}_\mathsf{M}(q)$ and thus $\chi_\mathsf{M}(q)$ do form a log-concave sequence. 
Of course we could have seen this for our example right away as $\chi_{\mathsf{K}_n}(q)$ has only real roots which implies log-concavity.

\subsection*{Acknowledgments}
\noindent
DC is supported by the SFB-TRR 195-- 286237555 ``Symbolic Tools in Mathematics and their Application''.
LK is supported by the SFB-TRR 358 -- 491392403 ``Integral Structures in Geometry and Representation Theory''.
BS is supported by the Swedish Research Council grant 2022-04224.

\bibliographystyle{siam}
\bibliography{references}
\end{document}